\documentclass[12pt]{article}
\usepackage[francais]{babel}
\usepackage{amsmath,amsfonts,amssymb,amsthm}
\newcommand{\Z}{\mathbb{Z}}
\newcommand{\F}{\mathcal{F}}
\newcommand{\T}{\mathbb{T}}
\newcommand{\N}{\mathbb{N}}
\newcommand{\Lgr}{\Longrightarrow}
\newcommand{\lgr}{\longrightarrow}
\newcommand{\rg}{\rightarrow}

\begin{document}

\title{Ensembles quasi--ind\'ependants\\
 et ensembles de Sidon, th\'eorie de Bourgain (d'apr\`es Myriam D\'echamps, Li et Queffelec \cite{lique},)}
\author{Jean--Pierre kahane\\
Laboratoire de Math\'ematique,\\
Universit\'e Paris--Sud \`a Orsay}
\date{}
\maketitle

Les ensembles de Sidon ont 50 ans. Leur pr\'ehistoire comprend les s\'eries lacunaires de Weierstrass, de Poincar\'e et d'Hadamard, utilis\'ees pour mettre en \'evidence des fonctions continues nulle part d\'erivables, et des  s\'eries de Taylor nulle part prolongeables en dehors du disque de convergence. La mise en \'evidence de leurs propri\'et\'es fondamentales dans le cas des suites lacunaires \`a la Hadamard est due \`a Sidon et \`a Banach. Leur d\'efinition appara\^{i}t en 1957 et 1960 dans deux articles, de Walter Rudin et de moi, et elle est popularis\'ee par le livre de Rudin ``Fourier analysis on groups''  \cite{kaha},  \cite{rudin},  \cite{rudin2}.

Les ensembles de Sidon ont servi de banc d'essai \`a la m\'ethode de s\'election al\'eatoire, introduite dans une note de Katznelson et Malliavin relative \`a la ``conjecture de dichotomie''  \cite{katmal}, et reprise par Katznelson en vue de l'\'etude du comportement des ensembles de Sidon plong\'es dans le groupe de Bohr  \cite{katz}. La s\'election al\'eatoire est essentielle dans la th\'eorie de Bourgain que je vais exposer.

La r\'eunion de deux ensembles de Sidon est un ensemble de Sidon : c'est le th\'eor\`eme de Drury (1970)  \cite{dru}, qui a \'et\'e le point de d\'epart de travaux brillants jusqu'au milieu des ann\'ees 80. Le recours aux probabilit\'es appara\^{i}t  d\`es la note de Drury et il est pr\'esent dans tous les travaux ult\'erieurs.

Il y a une parent\'e manifeste entre les ensembles de Sidon et les suites de variables al\'eatoires ind\'ependantes et \'equidistribu\'ees. Mais dans la plupart des groupes ab\'eliens discrets (par exemple $\Z$) il n'existe pas de sous--ensemble ind\'ependant infini. Un outil de base pour la mise en  \'evidence et l'\'etude des ensembles de Sidon est constitu\'e par les produits de Riesz, comme on va le voir.  Et le cadre naturel pour les produits de Riesz est celui des ensembles quasi--ind\'ependants, qui sont des ensembles de Sidon particuliers.

Toute r\'eunion finie d'ensembles quasi--ind\'ependants est un ensemble de Sidon ; c'est une cons\'equence du th\'eor\`eme de Drury, et la preuve directe est facile (voir ci--apr\`es). Est--il vrai que tout ensemble de Sidon est une r\'eunion finie d'ensembles quasi--ind\'ependants ? La question est toujours ouverte.

Un th\'eor\`eme de Pisier caract\'erise les ensembles de Sidon en termes d'ensembles quasi--ind\'ependants. Le but principal de cet article est d'en faire l'expos\'e, la g\'en\'eralisation et la d\'emonstration par la m\'ethode de Bourgain. C'est le parti--pris de Daniel  Li et Herv\'e Queffelec dans leur livre, et ils se r\'ef\`erent \`a Myriam D\'echamps pour expliciter la pens\'ee de Jean Bourgain, exprim\'ee dans ses articles de fa\c con tr\`es concentr\'ee. Le pr\'esent article n'apporte essentiellement rien de nouveau. Sa r\'edaction a suivi une invitation de Ben Green \`a l'Institut Newton de Cambridge, \`a tenter de faire le point sur les ensembles de Sidon, en janvier 2007.

Les r\'ef\'erences principales se trouvent dans  \cite{lique}. Je me borne ici \`a  \cite{pisier},  \cite{pisier2},  \cite{pisier3} ,  \cite{bourg}, pour les travaux de Pisier et Bourgain.

La partie 1 introduira les notations et d\'efinitions. La partie~2 donnera les premiers r\'esultats sur ensembles quasi--ind\'ependants et ensembles de Sidon et la partie~3 quelques exemples. La partie~4, la plus importante, sera un expos\'e de la th\'eorie de Bourgain.

\section{Notations et premi\`eres d\'efinitions \'equivalentes}

$G$ est un groupe ab\'elien compact.

$\Gamma$ est son dual, groupe ab\'elien discret, not\'e multiplicativement.

$\Lambda$ est une partie de $\Gamma$.

$\Lambda$ est quasi--ind\'ependante $(qi)$ si l'\'egalit\'e $\prod \lambda^{\varepsilon(\lambda)}=1$, o\`u $\lambda\in \Lambda$, $\varepsilon(\lambda) \in \{-1,0,1\}$ et $\sum |\varepsilon(\lambda)|<\infty$, entra\^{i}ne que tous les $\varepsilon(\lambda)$ sont nuls. Dans toute la suite, on supposera sans le pr\'eciser $\lambda\in \Lambda$.

$\Lambda$ est $S$--Sidon si, pour toute somme finie $f=\sum \hat f (\lambda)\lambda$, on a
$$
\sum |\hat f(\lambda)| \le S \sup_{g\in G} |f(g)|
$$
ce que nous \'ecrirons sous la forme
$$
\| f\|_A \le S\|f\|_C\,;
$$
$A=A(G)$ est l'alg\`ebre de Banach constitu\'ee par les sommes de s\'eries de Fourier absolument convergentes, c'est--\`a--dire $A(G) =\F \ell^1(\Gamma)$ ; 

\noindent $C=C(G)$ est l'alg\`ebre de Banach constitu\'ees par les fonctions continues sur $G$ ; les normes sont explicit\'ees ci--dessus.

$\Lambda$ est un ensemble de Sidon (en bref, est Sidon) s'il est $S$--Sidon pour un certain r\'eel $S$. Sa constante de Sidon est la borne inf\'erieure des $S$ en question.

Les sous--espaces ferm\'es engendr\'es par $\Lambda$ dans $A$, $C$ et $L^\infty$ $(=L^\infty(G))$, sont not\'es $A_\Lambda$, $C_\Lambda$ et $L_\Lambda^\infty$.

``$\Lambda$ est Sidon''  s'exprime, de fa\c con \'equivalente, par l'une des \'egalit\'es d'ensembles qui suivent :
$$
\begin{array}{ll}
 A_\Lambda &  = C_\Lambda\,,   \\
A_\Lambda  & = L_\Lambda^\infty\,,   \\
  \ell^\infty(\Lambda)   &   = \widehat M\mid_\Lambda  \\
c_0(\Lambda) &=\widehat L^1\mid _\Lambda  
\end{array}
$$
dans lesquelles $\ell^\infty(\Lambda)$ et $c_0(\Lambda)$ repr\'esentent respectivement les suites (ou fonctions) born\'ees resp. tendant vers $0$ \`a l'infini sur $\Lambda$, et $\widehat M \mid_\Lambda$ resp.~$\widehat{L^1}\mid_\Lambda$ les restrictions \`a $\Lambda$ des transform\'ees de Fourier de mesures de Radon born\'ees resp. de fonctions int\'egrables sur~$G$.

Voici un crit\`ere commode : $\Lambda$ est Sidon sii (si et seulement si)
$$
\exists \delta >0 : \forall \varphi \in \{-1,1\}^\Lambda \ \exists \mu\in M(G) : \quad \forall\xi \lambda |\widehat \mu (\lambda) -\varphi (\lambda)| \le 1-\delta \,.
$$

Enfin, $\Lambda$ est $S$--Sidon sii pour toute application $\lambda\lgr z(\lambda)$, $|z(\lambda)|=1$, il existe $\mu\in M(G)$, $\int | d\mu|=1$, telle que $R e (z(\lambda) \int \lambda d\mu) \ge {1\over S}$.

Toutes ces \'equivalences sont faciles \`a \'etablir et se trouvent dans  \cite{rudin2}  et~ \cite{lique}.

\section{Quasi--ind\'ependants et ensembles de Sidon, premiers r\'esultats}

i) Tout $q i$ est Sidon

ii) Toute r\'eunion finie de $qi$ est Sidon

iii) Soit $\Lambda\subset \Gamma$. S'il existe $c=c(\Lambda)>0$ tel que, pour toute mesure positive $\varpi$ sur $\Lambda$, il existe une partie $\Lambda'$ de $\Lambda$ $qi$ et telle que $\varpi(\Lambda') \ge \varpi(\Lambda)$, $\Lambda$ est Sidon (et $S$--Sidon avec $S=S(c))$.

Voici les preuves, et quelques compl\'ements. Pour justifier les calculs, on peut se restreindre d'abord au cas o\`u $\Lambda$ est fini.

i) Supposons $\Lambda$ $qi$. Le produit de Riesz associ\'e est
$$
R= \prod (1+{1\over2} (\lambda +\overline\lambda)) = \prod (1+{1\over2} (\lambda+\lambda^{-1}))\,.
$$
Son d\'eveloppement est 
$$
R= \sum_{(\varepsilon_\lambda)} \Big({1\over2}\Big)^{\sum|\varepsilon_\lambda|}\prod \lambda^{\varepsilon_\lambda}\,, \qquad \varepsilon_\lambda = \{-1,0,1\}\,, \sum|\varepsilon(\lambda)| <\infty\,.
$$
Comme, suite \`a la quasi--ind\'ependance, le terme constant est 1, $R$ est une mesure de probabilit\'e et on peut \'ecrire $R=\sum \widehat R(\gamma)\gamma$ (le signe $=$ signifiant le d\'eveloppement formel en s\'erie de Fourier) avec $\widehat R(1)=1$, $0\le \widehat R(\gamma)\le 1$,
$$
\widehat R(\lambda) = {1\over2} +\widehat r_\lambda(\lambda)\qquad \hbox{avec }\ 0\le \widehat r_\lambda(\lambda)\le {1\over2}\,.
$$

Posons maintenant
$$
R_{a,z} = \prod \Big(1+{a\over2} (z_\lambda \lambda+ \overline z_\lambda\overline\lambda)\Big)
$$
o\`u $0 < a <1$ et $|z_\lambda|=1$ ;  $z=(z_\lambda)_{\lambda\in \Lambda}$. On a de nouveau une mesure de probabilit\'e, avec $\widehat R_{a,z}(1)=1$, et
$$
\widehat R_{a,z}(\lambda) = {a\over2} z_\lambda +\widehat r_{a,z,\lambda}(\lambda) \quad\hbox{avec}\quad |\widehat r_{a,z,\lambda}(\lambda)| \le {a^2\over2}\ \hbox{si}\quad\overline\lambda\neq \lambda\,.
$$
La consid\'eration de $R_{a,z}$ permet de conclure, mais on a une preuve plus nette et une meilleure constante de Sidon en introduisant
\[
\begin{array}{cl}
  R_{a,z}^*&= \displaystyle\int_0^{2\pi} R_{a,ze^{it}} e^{-it}{dt\over2\pi}  \\
  &=\displaystyle \sum_{(\varepsilon_\lambda)} \Big({a\over2}\Big)^{\sum|\varepsilon_\lambda|} \prod(z_\lambda \lambda)^{\varepsilon_\lambda} \int_0^{2\pi} \prod e^{it\varepsilon_\lambda} e^{-it} {dt\over2\pi}   \\
  &=\displaystyle\sum_{\Sigma \varepsilon_\lambda=1} \Big({1\over2}\Big)^{\sum|\varepsilon_\lambda|} \prod (z_\lambda \lambda)^{\varepsilon_\lambda}   
\end{array}
\]
C'est une mesure de masse totale $\le 1$, et
$$
\widehat R_{a,z}^* (\lambda) = {a\over2} z_\lambda +\widehat r_{a,z,\lambda}^*(\lambda)
$$
avec
$$
|\widehat r_{a,z,\lambda}^* (\lambda) | \le \widehat r_{a,1,\lambda}^* (\lambda) \le {a^3\over2}
$$
donc
$$
Re  (\overline z_\lambda R_{a,z}^* (\lambda)) \ge {a\over2} - {a^3\over2}\,.
$$
En choisissant $a={1\over\sqrt{3}}$, on voit d'apr\`es le dernier crit\`ere de la partie~1 que $\Lambda$ est $S$--Sidon avec $S=3\sqrt{3}$ $(<5,2)$.\hfill $\blacksquare$

On peut pousser l'\'etude, en regardant le cas $\sum|\varepsilon_\lambda|=3$, et on arrive~\`a
$$
Re (\overline z_\lambda R_{a,z}^* (\lambda)) \ge {1\over2} - {a^3\over8} - {a^5\over2}\,;
$$
et, en choisissant $a$ pour que le second membre soit minimum, on voit que $\Lambda$ est $S$--Sidon avec $S=4,27$.

La valeur optimale de $S$ n'est pas connue.

Remarquons que, pour $\gamma\notin \Lambda$,
$$
|\widehat R_{a,z}^* (\gamma) | \le \widehat R_{a,1}^* (\gamma) \le a^3 \widehat R(\gamma)\le a^3\,.
$$

ii) Si $\Lambda$ est une r\'eunion de $k$ ensembles $qi$, on peut les supposer disjoints, soit $\Lambda',\Lambda'',\ldots \Lambda^{(k)}$, et, avec des notations \'evidentes on pose
$$
R_{a,z}^{**} = {1\over k} (R_{a,z'}^* +R_{a,z''}^* + \cdots R_{a,z^{(k)}}^*)
$$
d'o\`u
$$
Re (\overline z_\lambda \widehat R_{a,z}^{**} (\lambda)) \ge {1\over k} \Big({a\over2} - {a^3\over2} - (k-1) a^3\Big)
$$
et, en choisissant $a$ de fa\c con \`a minimiser le second membre, $\Lambda$ est $S$--Sidon avec
$$
\hbox to 5cm{} S= 3\sqrt{3}\ k\ \sqrt{2k-1}\,.\hbox to 4,5cm{} \blacksquare
$$

iii) Supposons v\'erifi\'ee la condition de Bourgain

\textbf{(CB)} : il existe $c=c(\Lambda)>0$ tel que, pour toute mesure positive $\varpi$ sur $\Lambda$, il existe une partie $\Lambda'\subset \Lambda$, $qi$ et telle que $\varpi(\Lambda') \ge c\varpi(\Lambda)$.

Soit $f=\sum \widehat f (\lambda)\lambda$, $\sum|\widehat f (\lambda)|<\infty$. Choisissons $\varpi(\lambda) = |\widehat f(\overline\lambda)|$, puis, selon \textbf{(CB)}, $\Lambda'\subset \Lambda$, $qi$, avec
$$
\sum_{\lambda\in \Lambda'} |\widehat f(\overline\lambda)| \ge c \sum_{\lambda\in \Lambda} |\widehat f(\overline\lambda)|\,.
$$
On construit $R_{a,z}^*$ ci--dessus \underline{sur $\Lambda'$}, en prenant  $\overline z_\lambda \widehat f(\lambda) =|\widehat f(\lambda)|$ pour $\lambda\in \Lambda'$. Alors
\[
\begin{array}{rl}
\|f\|_C  &\displaystyle\ge Re \int f R_{a,z}^* = \sum_{\lambda\in\Lambda'} |\widehat f(\overline\lambda)| \ Re(\bar z_\lambda R_{a,z}^*(\lambda)) + \sum_{\lambda\in \Lambda\backslash \Lambda'} \hbox{idem} \\
  &\displaystyle\ge \sum_{\lambda\in\Lambda'} |\widehat f(\overline \lambda)| \Big({a\over2}-{a^3\over2}\Big) - \sum_{\lambda\in \Lambda}|\widehat f(\overline\lambda)|a^3  \\
  &\displaystyle\ge \Big(c\Big({a\over2}- {a^3\over2} \Big) -a^3\Big) \|f\| _A\,,   
\end{array}
\]
soit, en prenant la plus petite valeur de la parenth\`ese, $\|f\|_C\ge {1\over S}\|f\|_A$ avec 
$$
\hbox to 4cm{}{1\over S} = {1\over12\sqrt{3}}\ c^{3/2} (c+2)^{-1/2}(6-c) \hbox to 3cm{}\blacksquare
$$

\section{Exemples et remarques}

Prenons $G=\T$. Alors $\Gamma=\Z$, qu'on a coutume de noter additivement. L'ensemble $\{2^j, j\in \N\}$ est $q.i$. Il en est de m\^eme pour $\{\lambda_j\}$ si les $\lambda_j$ sont des entiers tels que $\lambda_{j+1}/\lambda_j \ge 2$ $(j\in \N)$. $\Lambda=\{\lambda_j\}$ est ``lacunaire \`a la Hadamard'' si $\lambda_{j+1}/\lambda_j \ge q$ $(j\in \N)$ pour un $q>1$ ; c'est alors une r\'eunion finie d'ensembles du type pr\'ec\'edent, donc c'est un Sidon.

A ma connaissance, on ne conna\^{i}t  explicitement la constante de Sidon d'un ensemble d'entiers (lorsqu'elle est finie) que dans quelques cas d'ensembles finis  \cite{neu}.

Rappelons la notation usuelle pour $f=\sum \widehat f (\gamma)\gamma$ :
$$
f(t) = \sum_{n\in \Z} \widehat f(n) \ e^{2\pi int} \qquad (t \in \T)\,.
$$

Prenons maintenant $G=\T^\N$, $\Gamma=\Z^{\oplus \N}$ (partie de $\Z^\N$ constitu\'e des $(n_j)_{j\in \N}$ tels que $\sum |n_j| <\infty$). Souvent $G$ est pris comme espace de probabilit\'e, $\Omega$, et ses \'el\'ements sont alors not\'es $(\omega_0,\omega_1,\ldots \omega_j,\ldots)$. Les $e^{2\pi i\omega_j}$ $(j\in\N)$ sont des variables al\'eatoires ind\'ependantes ; on les appelle $v\cdot a\cdot$ de Steinhaus et les s\'eries $\displaystyle\sum_{j\in \N}a_j e^{2\pi i\omega_j}$ s\'eries de Steinhaus. Elles sont de la forme $\displaystyle \sum_{\lambda\in \Lambda} \widehat f(\lambda)\lambda$ en prenant pour $\Lambda$ l'ensemble des vecteurs de base de $\Z^{\oplus\N} : (1,0,0,\ldots)$, $(0,1,0,\dots)$, $(0,0,1,\ldots)\ldots,$ ou plus simplement en notation multiplicative l'ensemble des $v\cdot a\cdot$ de Steinhaus. Il est imm\'ediat que $\Lambda$ est $qi$ et que sa constante de Sidon est $S=1$. 

Prenons enfin $G=(\Z/2\Z)^\N$. Les caract\`eres sur  $G$ sont les fonctions de Walsh, \`a valeurs $\pm1$, engendr\'ees par les fonctions de Rademacher $r_j=(-1)^{g_j}$ $(g=(g_0,g_1,\ldots))$. Les s\'eries de Rademacher sont les s\'eries $\displaystyle\sum_{j\in \N} a_jr_j$ $(=\sum \pm a_j)$. L'ensemble $\Lambda$ des fonctions de Rademacher, qui sont ind\'ependantes, est un ensemble $qi$. On conna\^{i}t sa constante de Sidon, qui est $\pi\over2$~ \cite{sei}.

Voici un r\'esum\'e de la preuve.

Par d\'efinition, $S$ est la borne sup\'erieure des $\sum |a_j|$ pour tous les choix de $a_j$ complexes nuls \`a partir d'un certain rang, tels que toutes les sommes $\sum \pm a_j$ soient de module $\le 1$. Nous pouvons nous borner, quitte \`a r\'eajuster les $a_j$, au cas
$$
\sup_{(\pm)} |\sum \pm a_j| = \sum a_j =1\,.
$$
La premi\`ere \'egalit\'e implique $Re a_j>0$ soit $-{\pi\over2}\le \arg a_j <{\pi\over2}$ pour tout $j$. Ordonnons les $a_j$ par arguments d\'ecroissants ; ainsi les points $0,a_1, a_1+a_2,\ldots a_1+a_2,\ldots,1$ sont  les sommets cons\'ecutifs d'une ligne polygonale concave. Tous ces points sont dans le demi--disque sup\'erieur de diam\`etre $[0,1]$ ; en effet, si l'un d'eux, $a_1+a_2+\cdots a_k$, \'etait \`a l'ext\'erieur, il verrait le segment $[0,1]$ sous un angle $<{\pi\over2}$, et il s'ensuivrait que $|a_1+a_2 +\cdots +a_k -(a_{k+1}+a_{k+2}+\cdots)|>1$. La  longueur de la ligne polygonale est $\sum |a_j|$. On a donc $\sum |a_j| <{\pi\over2}$, et en choisissant une ligne proche du demi--cercle, on voit que $\sup\limits_{(a_j)} \sum |a_j| = {\pi\over2}$. 

\section{Les conditions de Rudin, Pisier et Bourgain}

Ce sont  trois conditions n\'ecessaires et suffisantes pour que $\Lambda$ soit Sidon. Nous les d\'esignons par \textbf{(CR)},\textbf{ (CP)} et \textbf{(CB)}.

\vskip2mm

\begin{itemize}
\item[\textbf{(CR)}] \textit{Il existe un $C>0$ tel que, pour tout polyn\^ome} $f=\displaystyle\sum_{\lambda\in \Lambda}\widehat f(\lambda)\lambda$,
$$
\|f\|_p = C \sqrt{p} \|f\|_2 \qquad \hbox{pour }\ p>2
$$
\item[\textbf{(CP)}] \textit{Il existe un $b>0$ tel que toute partie finie $A$ de $\Lambda$ contienne un ensemble $qi$ $B$ tel que $|B| \ge b|A|$}
\item[\textbf{(CB)}] \textit{(d\'ej\`a \'ecrite en 2 iii). Il existe $c>0$ tel que, pour toute mesure positive $\varpi$ sur $\Lambda$, il existe une partie $\Lambda'$ de $\Lambda$, $qi$ et telle que $\varpi (\Lambda') \ge c\varpi(\Lambda)$.
}\end{itemize}

\vskip2mm

\textbf{\textit{Historique et remarques}}

\textbf{(CR)} a \'et\'e \'etabli par Rudin (1960  \cite{rudin}) comme condition n\'ecessaire, et par Pisier  \cite{pisier} comme condition suffisante. \textbf{(CP)} a \'et\'e \'etabli par Pisier comme condition n\'ecessaire et suffisante  \cite{pisier2}, \cite{pisier3}.  \textbf{(CB)}, qui implique \textbf{(CP)} en prenant pour $\varpi$ la mesure de d\'ecompte, a \'et\'e introduit par Bourgain \cite{bourg}, ainsi que l'encha\^{i}nement
$$
\hbox{Sidon} \Longrightarrow \hbox{\textbf{(CR)}} \Longrightarrow \hbox{\textbf{(CP)}} \Longrightarrow  \hbox{\textbf{(CB)}} \Longrightarrow  \hbox{Sidon}\,. 
$$
Nous avons d\'ej\`a vu en 2 \textbf{(CB)} $ \Longrightarrow$ Sidon. Reste \`a montrer Sidon $ \Longrightarrow$ \textbf{(CR)}, \textbf{(CR)} $ \Longrightarrow$ \textbf{(CP)} et \textbf{(CP) }$ \Longrightarrow$ \textbf{(CB)}.

\vskip4mm

\textbf{{\large Sidon implique (CR)}}

\vskip2mm

\textbf{\textsl{Proposition pr\'eliminaire sur les normes}} $L^q$ $(q\ge 1)$. \textit{Soit $\Lambda = \{\lambda_k\}$ un Sidon de constante $<S$, et $\Omega =\T^N = \{(\omega_k)\}$. Quels que soient les $a_k$ complexes, et $q\ge 1$,}
$$
{1\over S}\ \Big \|\sum a_k \ e^{2\pi i\omega_k} \Big\|_{L^q(\Omega)} \le \Big \|\sum a_k \lambda_k \Big\|_{L^q(G)}\le S\ \Big\|\sum a_k e^{2\pi i\omega_k}\Big\|_{L^q(\Omega)}\,.
$$

\vskip2mm

\textbf{\textit{Preuve}} : on introduit
\[
\begin{array}{ll}
  & \mu_\omega \in M(G),\ \|\mu_\omega\| \le S, \ \widehat{\mu_\omega} (\lambda_k) = e^{2\pi i\omega_k}   \\
  \noalign{\vskip1mm}
  &\widetilde{ \mu_\omega }\in M(G),\ \| \widetilde{\mu_\omega}\| \le S, \ \widehat{\widetilde{\mu_\omega}} (\lambda_k) = e^{2\pi i\omega_k}      \\
    \noalign{\vskip1mm}
  & f = \sum a_k \lambda_k\,,\ \ f_\omega = \sum a_k e^{2\pi i\omega_k}\lambda_k \,; 
\end{array}
\]
ainsi $f_\omega = f \ast \mu_\omega$ et $f= f_\omega \ast \widetilde{\mu_\omega}$, donc
$$
\|f\|_{L^q(G)} \le \|\widetilde{\mu_\omega}\|\ \|f_\omega\|_{L^q(G)} \le \|f_\omega\|_{L^q(G)}
$$
soit
$$
\int_G |f|^q \le S^q \int_G |f_\omega|^q\,,
$$
et, en prenant la moyenne sur $\Omega$,
$$
\int_G |f|^q \le S^q \ E \int_G |f_\omega|^q = S^q \int_G \ E |f_\omega|^q = S^q \ E \ \Big| \sum a_k e^{2\pi i\omega_k}\Big|^q
$$
et de m\^eme pour la premi\`ere in\'egalit\'e.

\vskip2mm

\textbf{\textit{Preuve que Sidon implique} (CR)} (in\'egalit\'es de Rudin) :

On \'etablit d'abord, ce qui est classique, et facile au moyen de la transform\'ee de Laplace, que pour des $a_k$ r\'eels
$$
\Big\|\sum \pm  a_k\Big\|_{L^q(\Omega)} \le c_0 \sqrt q \Big(\sum a_k^2\Big)^{1/2}
$$
puis on utilise deux fois la proposition pr\'eliminaire $((\pm) \rg (\omega_k)\rg (\lambda_k))$ pour conclure.

\vskip 2mm

\textbf{\large (CR) implique (CP)}

\vskip1mm

Pour simplifier les \'ecritures nous allons supposer $\lambda^2\neq 1$ pour tout $\lambda\in \Lambda$, en plus de l'hypoth\`ese \textbf{(CR)}. Nous \'ecrirons $\Lambda$ au lieu de la partie finie~$A$.

Si $D \subset \Lambda$, et si $d$ est un entier $\neq1$, le nombre de relations de hauteur $d$ dans $D$, c'est--\`a--dire d'\'egalit\'es $\prod\limits_{\lambda\in D}\lambda^{\varepsilon_\lambda}=1$,  avec $\varepsilon_\lambda\in \{-1,0,1\}$ et $\sum |\varepsilon_\lambda|=d$,~est
$$
\int_G \prod_{\lambda\in D, \alpha_\lambda\in\{0,1\},\sum \alpha_\lambda=d} \ \ (\lambda+\overline{\lambda})^{\alpha_\lambda}\,.
$$
Choisissons une suite de variables al\'eatoires de Bernoulli ind\'ependantes $\xi_\lambda$ $(\lambda\in \Lambda)$ de m\^eme loi $P(\xi_\lambda=1)=\eta$, $P(\xi_\lambda=0)=1-\eta$, et d\'efinissons $D$ par s\'election al\'eatoire :
$$
D=D(\omega) = \{\lambda\in \Lambda : \xi_\lambda=1\}
$$
Le nombre de relations de hauteur $d$ dans $D$ est alors 
$$
\int_G \prod_{\lambda\in \Lambda,\alpha_\lambda\in\{0,1\},\sum \alpha_\lambda=d} (\lambda+\overline{\lambda})^{\alpha_\lambda\xi_\lambda}
$$
et son esp\'erance est 
$$
\eta^d \int_G  \prod_{\lambda\in \Lambda,\alpha_\lambda\in\{0,1\},\sum \alpha_\lambda=d} (\lambda+\overline{\lambda})^{\alpha_\lambda}\,.
$$
Or
$$
 \prod_{\lambda\in \Lambda,\alpha_\lambda\in\{0,1\},\sum \alpha_\lambda=d} (\lambda+\overline{\lambda})^{\alpha_\lambda} \le {1\over d!} \Big(\sum_{\lambda\in \Lambda} (\lambda+\overline{\lambda})\Big)^d\,.
$$
L'hypoth\`ese \textbf{(CR)} dit que
\[
\begin{array}{rl}
\displaystyle  \int_G \Big(\sum_\mu\lambda\in \Lambda (\lambda+\overline{\lambda})\Big)^d&\displaystyle \le C^d d^{d/2}\Big(\int_G \Big(\sum_{\lambda\in \Lambda}(\lambda+\overline{\lambda})\Big)^2\Big)^{d/2}   \\
  &= C^d\ d^{1/2} (2|\Lambda|)^{d/2}\,. \\
\end{array}
\]
L'esp\'erance du nombre de relations de hauteur $>\ell$ dans $D$ est major\'ee par 
$$
\sum_{d>\ell} \eta^d {1\over d!} C^d d^{d/2}(2|\Lambda|)^{d/2} \le \Big(2\eta\ C\ e \sqrt{|\Lambda|\over\ell} \Big)^\ell
$$
Elle est $<2^{-\ell}$ si $\eta= {1\over 4Ce}$ (nous faisons ce choix) et $\ell \ge {1\over4} \eta|\Lambda|$.

Choisissons $\ell = {1\over4} \eta|\Lambda|$. Comme $E|D| =\eta |\Lambda|$, il existe un $\omega$ tel que $|D| > {1\over2} \eta |\Lambda|$ et qu'il n'y ait dans $D$ aucune relation de hauteur $>\ell$ ; remarquons que $\ell < {1\over2}|D|$ pour le $D=D(\omega)$ choisi.

Alors vient la belle id\'ee. On prend dans $D$ une relation de hauteur maximale $(\le \ell)$. C'est du type $\prod\limits_{\lambda\in D'}\lambda^{\varepsilon_\lambda}=1$ avec $\varepsilon_\lambda=\pm1$, pour un $D'\subset D$, tel que $|D'|\le \ell$. Choisissons $B= D \setminus D'$. S'il y avait une relation dans $B$, on en ferait le produit membre \`a membre avec la relation de hauteur maximale, et comme $B$ et $D'$ sont disjoints on obtiendrait une relation de hauteur strictement sup\'erieure, ce qui est impossible. Donc $B$ est~$qi$.

De plus $|B| >{1\over2}|D| > {1\over4} \eta |\Lambda| = {1\over16Ce}|\Lambda|$, donc $|B| > b|\Lambda|$ $(=b|A|)$ avec $b={1\over16C e}$ (ind\'ependant de~$A$).

Le cas o\`u $\lambda^2=1$ pour certains $\lambda\in\Lambda$ ne n\'ecessite que des modifications d'\'ecriture que je laisse au lecteur. On a donc montr\'e \textbf{(CR)} $\Lgr$ \textbf{(CP)}. 

\vskip2mm

\textbf{\large (CP) implique (CB)}

\vskip1mm

C'est la partie la plus laborieuse de la th\'eorie de Bourgain. Je la pr\'esenterai en trois \'etapes : 1) r\'eduction de la mesure $\varpi$ \`a une forme plus maniable, \`a savoir une combinaison lin\'eaire de mesures de d\'ecompte sur des parties $\Lambda_j$ de $\Lambda$, $qi$ et de tailles tr\`es diff\'erentes 2) construction par s\'election al\'eatoire de parties $A_j$ des $\Lambda_j$ dont la d\'ependance \`a l'\'egard des autres $\Lambda_k$ est soigneusement contr\^ol\'ee 3) utilisation d'un argument de maximalit\'e analogue \`a celui de la preuve de \textbf{(CR)}~$\Lgr$~\textbf{(CP)} pour obtenir des parties $\Lambda_j'$ des $A_j$ dont la r\'eunion est $\Lambda'$, l'ensemble $qi$ cherch\'e. Le jeu consistera \`a conserver une proportion notable de la masse \`a chaque \'etape.

%\vskip2mm

\eject
\textbf{\textit{Premi\`ere \'etape.}}\ Elle se d\'eroule en plusieurs temps.

\noindent 1.1. \hskip2mm Posons
$$
\varpi_1(\lambda) = \sup_{k\in \Z}\  \{2^{-k} : 2^{-k} \le \varpi (\lambda)\}\,.
$$
Ainsi $\varpi_1\Lambda) \ge {1\over2} \varpi (\Lambda)$. Soit $A_k = \{\lambda : \varpi_1(\lambda) = 2^{-k}\}$ ; alors
$$
\varpi_1 = \sum_{k\in\Z} 2^{-k} \sum_{\lambda\in A_k} \delta _\lambda
$$
(la derni\`ere somme est la mesure de d\'ecompte sur $A_k$).

\vskip2mm

\noindent1.2. \hskip2mm Utilisant \textbf{(CP)}, soit $B_k \subset A_k$, $qi$, avec $|B_k| \geq b|A_k|$. Posons
$$
\varpi_2 = \sum_{k\in \Z} 2^{-k} \sum_{\lambda\in B_k} \delta _\lambda
\,.$$
Ainsi $\varpi_2(\Lambda) \ge b \ \varpi_1(\Lambda)$.

\vskip2mm

\noindent1.3.  \hskip2mm Soit $R>1$ ($R$ sera d\'efini \`a la seconde \'etape) et
$$
K_j = \{k : R^j \le |B_k| < R^{j+1}\} \qquad (j=0,1,2,\ldots)\,.
$$
Posons $k(j) =\inf K_j$ $(=-\infty$ si $K_j=\emptyset)$, et
$$
\varpi_3 = \sum_j 2^{-k(j)} \sum_{\lambda\in B_{k(j)}} \delta _\lambda\,.
$$
Si $k\in K_j$,
$$
\varpi_2 (B_k) = 2^{-k} |B_k| \le \varpi_3(B_{k(j)}) 2^{-(k-k(j))}R
$$
donc
$$
\varpi_2 \Big(\bigcup_{k\in K_j}B_k \Big) \le 2\ R\ \varpi_3 (B_{k(j)})
$$
et $\varpi_3(\Lambda) \ge {1\over 2R} \varpi_2(\Lambda) $.

\vskip2mm

\noindent 1.4.\hskip2mm Si
$$
\varpi_3  \Big(\bigcup_{j\ \hbox{pair}}B_{k(j)} \Big) \ge \varpi_3  \Big(\bigcup_{\hbox{impair}}B_{k(j)} \Big) 
$$
on pose $ \Lambda_j = B_{k(2j)}$ $(j=0,1,2,\ldots)$ et
$$
\varpi_4 = \sum_j 2^{-k(2j)} \sum_{\lambda\in \Lambda_j} \delta _\lambda\,.
$$
Sinon, on pose $ \Lambda_j=B_{k(2j-1)}$ $(j=1,2,\ldots)$ et
$$
\varpi_4 =\sum_j 2^{-k(2j-1)} \sum_{\lambda\in  \Lambda_j} \delta _\lambda\,.
$$
Dans les deux cas, $\varpi_4( \Lambda) \ge {1\over2}\varpi_3( \Lambda)$.

\vskip2mm

\noindent 1.5. \hskip2mm Si $\varpi_4( \Lambda_0) \ge {1\over2} \varpi_4( \Lambda)$ on choisit $\Lambda'=\Lambda_0$ et on a \'etabli \textbf{(CB)}, avec $c={b\over16R}$.  On va donc se restreindre dans la suite au cas $\varpi_4( \Lambda_0) <{1\over2}\varpi_4( \Lambda)$, et on d\'efinit $\varpi_5$ comme la restriction de $\varpi_4$ \`a $\Lambda\setminus  \Lambda_0$ ; ainsi $\varpi_5(\Lambda) \ge {1\over2} \varpi_4(\Lambda) \ge {b\over16R}\varpi(\Lambda)$, $\varpi_5$ est une combinaison lin\'eaire \`a coefficients positifs de mesures de d\'ecompte sur les $\Lambda_j$, on a $|\Lambda_j| \ge R$ et ${|\Lambda_{j+1}\over |\Lambda_j|} \ge R$ pour $j=1,2,\ldots$.

Pour d\'emontrer \textbf{(CB)}, il suffira de trouver des $\Lambda_j'\subset \Lambda_j$, tels que $|\Lambda_j'| \ge {1\over10}|\Lambda_j|$, et que $\cup \Lambda_j'$ soit~$qi$.

\vskip2mm

\textbf{\textit{Seconde \'etape.}}\  Rappelons que les $\Lambda_j$ sont $qi$, avec $|\Lambda_1| \ge R$ et ${|\Lambda_{j+1}|\over|\Lambda_j|} \ge  R$ $(j=1,2,\ldots)$.

Fixons $j$. Soit $(\xi_\lambda)$ $(\lambda\in \Lambda_j)$ une suite de $v\cdot a\cdot$ de Bernoulli ind\'ependantes, avec $P(\xi_\lambda=1)={1\over4}$ et $P(\xi=0)={3\over4}$, et soit
$$
A(\omega) = \{\lambda\in \Lambda_j : \xi_\lambda =1\}\,.
$$
Comme $|\Lambda_j| \ge R$, on a $|A(\omega)| > {1\over5}|\Lambda_j|$ avec une probabilit\'e $p_R$ voisine de $1$ quand $R$ est~grand.

Soit $\rho$ un \'el\'ement de $\Gamma$ engendr\'e par $\bigcup\limits_{k\neq j}\Lambda_k$, c'est--\`a--dire
\[
\begin{array}{rl}
  \rho& =\displaystyle \prod_{k\neq j} \rho_k \qquad \hbox{(produit fini)}     \\
\rho_k  &= \displaystyle\prod_\ell \lambda_{k\ell}^{n_{k\ell}} \qquad (\lambda_{k\ell} \in \Lambda_k,\ n_{k\ell}\in \Z)\,.   \\
\end{array}
\]
Posons $d(\rho_k) = \sum_\ell |n_{k\ell}|$. Pour le moment, $\rho$ est fix\'e.

On va s'int\'eresser aux $\sigma$ engendr\'es par $A(\omega)$, de la forme $\sigma= \prod\limits_{\lambda\in A(\omega)} \lambda^{\varepsilon_\lambda}$ avec $\varepsilon_\lambda \in \{-1,0,1\}$, tels que $d(\sigma) = \sum |\varepsilon_\lambda| > L$ fix\'e. D\'esignons par $N(\omega)$ le nombre total de relations $\sigma\rho=1$. Ainsi
$$
N_\rho(\omega) = \int_G S_L \prod_{\lambda\in A(\omega)} (1+\lambda+\overline{\lambda})\rho
$$
$S_L\prod$ d\'esignant la somme des termes du d\'eveloppement dont la hauteur d\'epasse $L$. On peut \'ecrire
\[
\begin{array}{rl}
  N_\rho(\omega)&={\displaystyle\int_G S_L \prod_{\lambda\in \Lambda_j} (1+(\lambda+\overline{\lambda})\xi_\lambda)\rho}  \\
EN_\rho(\omega)  &={ \displaystyle\int_G S_L \prod_{\lambda\in \Lambda_j} \Big( 1+{1\over4}(\lambda+\overline{\lambda})\Big) \rho}\,.   \\
\end{array}
\]
La quasi--ind\'ependance de $\Lambda_j$ nous dit que $\prod\limits_{\lambda\in \Lambda_j}(1+{1\over2}(\lambda+\overline{\lambda}))$ est une mesure de probabilit\'e ; \'ecrivons son d\'eveloppement sous la forme $\sum ({1\over2})^{d(\sigma)}\sigma$. Alors
$$
EN_\rho (\omega) = \sum_{\sigma:d(\sigma)>L} \Big({1\over4}\Big)^{d(\sigma)} \int \rho\ \sigma \le {1\over2^L} \sum_\sigma \Big({1\over2}\Big)^{d(\sigma)} \int \rho\sigma \le {1\over2^\lambda}\,.
$$

Consid\'erons maintenant tous les $\rho=\prod\limits_{k\neq j}\rho_k$ pour lesquels $d(\rho_k)\le d_k$, suite positive fix\'ee, et soit $N(\omega)$ le nombre total des relations $\rho\sigma=1$ entre ces $\rho$ et les $\sigma$ consid\'er\'es ci--dessus. On a $N(\omega)=\sum_\rho N_\rho(\omega)$,
$$
EN(\omega) \le {1\over2^L}\ \sharp\ \{\rho\}\,.
$$
Or $\sharp\{\rho\}=\prod\limits_{k\neq j} \sharp \{\rho_k\}$. Evaluons $\sharp\{\rho_k\}$. C'est au plus le nombre de solutions de $\sum\limits_\ell |n_{k\ell}| \le d_k$ (voir ci--dessus la d\'efinition de $\rho_k$), que je d\'esigne par $N(d_k,|\Lambda_k|)$. Ainsi
$$
EN(\omega) \le {1\over2^L}\prod_{k\neq j} N(d_k,|\Lambda_k|)\,.
$$

Evaluons la fonction $N(d,q)$. C'est le nombre des fonctions $f:\{1,2,\ldots q\} \rg \Z$ telles que $\sum\limits_m |f(m)|\le d$. Comme le nombre de valeurs $f(m)\neq0$ ne d\'epasse pas $d\wedge q$, et que le nombre des fonctions $|f| :  \{1,2,\ldots q\}\rg \Z^+$ telles que $\sum |f|(m)\le d$ est exactement $
\left(
\begin{array}{cc}
  d+q   \\
  q   \\  
\end{array}
\right)
$ (on le voit en codant ces fonctions par des parties \`a $d$ \'el\'ements de $\{1,2,\ldots q+d\}$), on obtient
$$
N(d,q) \le 2^{d\wedge q}
\left(
\begin{array}{c}
  d+q   \\
  q   \\  
\end{array}
\right) \left( =2^{d\wedge q} 
\left(
\begin{array}{c}
  d+q   \\
  d   \\  
\end{array}
\right)\right)\, ,
$$
donc, pour une constante absolue $C$ ($C=20$ convient)
\[
\begin{array}{rll}
  N(d,q)&  \le C\Big({\displaystyle{q\over d}}\Big)^d & \hbox{si}\ d\le q  \\
 N(d,q)&  \le C\Big({\displaystyle{d\over q}}\Big)^q &  \hbox{si}\ q\le d\,. \\
\end{array}
\]
Ainsi
$$
EN(\omega) \le {1\over2^L} \exp \Bigg(\sum_{d_k\le |\Lambda_k|} d_k \log\Big(C{|\Lambda_k|\over d_k}\Big) + \sum_{d_k>|\Lambda_k|} |\Lambda_k| \log\Big(C{d_k\over|\Lambda_k|}\Big)\Bigg)\,.
$$

Choisissons $d_k|\Lambda_k| = |\Lambda_j|^2$ pour tout $k\neq j$. Ainsi $d_k\le |\Lambda_k|$ signifie $|\Lambda_j|^2 \le |\Lambda_k|^2$, c'est--\`a--dire $j<k$, donc
$$
EN(\omega) \le {1\over2^L} \exp \Bigg(\sum_{k\ge j} {|\Lambda_j|^2\over|\Lambda_k|} \log \Big(C{|\Lambda_k|^2\over|\Lambda_j|^2}\Big) +\sum_{k<j} |\Lambda_k| \log\Big(C {|\Lambda_j|^2\over |\Lambda_k|^2}\Big)\Bigg)\,.
$$
On peut alors faire usage de l'hypoth\`ese sur les quotients de $|\Lambda_j|$, et on obtient pour $R>R(C)$
$$
EN(\omega) \le {1\over2^L} \exp |\Lambda_j| \Big( \sum_{k\ge j} {1\over R^{k-j}} \log (CR^{2(k-j)}) + \sum_{k<j} {1\over R^{j-k}} \log (CR^{2(j-k)})\Big)\,.
$$

La parenth\`ese est aussi petite qu'on veut par choix de $R$. Choisissons $L={|\Lambda_j|\over 10}$. Alors, par choix de $R$, $EN(\omega)$ est aussi petit qu'on veut, donc \'egalement $P(N(\omega)\neq 0)$.

Donc, \`a condition de prendre $R$ assez grand, on peut choisir $\omega$ de fa\c con \`a avoir simultan\'ement $|A(\omega)| > {1\over5}|\Lambda_j|$ et $N(\omega)=0$, $N(\omega)$ d\'esignant le nombre total de relation $\sigma\rho=1$ entre les $\sigma=\prod\limits_{\lambda\in A(\omega)} \lambda^{\varepsilon_\lambda}$ $(\varepsilon_\lambda \in \{-1,0,1\})$ de hauteurs $\ge {|\Lambda_j|\over10}$ et les $\rho$ produits de $\rho_k$ engendr\'es par les $\Lambda_k$ $(k\neq j)$, de hauteurs $\le d_k = {|\Lambda_j|^2\over|\Lambda_k|}$.

D\'esormais, nous d\'esignons par $A_j$ cet $A(\omega)$ ainsi choisi.

\vskip2mm

\textbf{\textit{Troisi\`eme \'etape.}}\  Il s'agit de r\'ealiser le programme trac\'e \`a la fin de la premi\`ere \'etape. Commen\c cons par d\'efinir $\Lambda_j'$. Pour cela, fixons $j$, et consid\'erons tous les $\sigma=\prod\limits_{\lambda\in A_j} \lambda^{\varepsilon_\lambda}$ $(\varepsilon_\lambda \in\{-1,0,1\}$) et $\rho_k = \prod\limits_{\lambda\in \Lambda_k}\lambda^{n_\lambda}$ $(n_\lambda\in \Z)$ $(k\neq j)$, en nombre fini,  tels que $d(\rho_k)|\Lambda_k| \le d(\sigma)|\Lambda_j|$ $(k\neq j)$ et $\sigma\prod{\rho_k}=1$. Nous dirons alors que $(\sigma,(\rho_k))$ est un syst\`eme permis.

Comme $d(\sigma) |\Lambda_j|\le |\Lambda_j|^2$, nous savons par la conclusion de la seconde \'etape que les $\sigma$ d'un syst\`eme permis v\'erifient $d(\sigma) < {|\Lambda_j|\over10} < {1\over2}|A_j|$.

Distinguons deux cas. $\alpha)$ S'il n'existe pas de syt\`eme permis, posons $\Lambda_j'=A_j$. $\beta)$ S'il en existe, choisissons--en un tel que $d(\sigma)$ soit maximum, fixons $\sigma$ et $S_j$ son ``support'', c'est--\`a--dire l'ensemble minimal dans $A_j$ qui l'engendre : $\sigma=\prod\limits_{\lambda\in S_j} \lambda^{\alpha_\lambda}$ $(\alpha_\lambda = \pm 1)$. Posons alors $\Lambda_j'=A_j\setminus S_j$. On a bien dans tous les cas $|\Lambda_j'| >{1\over2} |A_j| >{1\over10} |\Lambda_j|$.

Reste \`a montrer que $\bigcup \Lambda_i'$ est $qi$. Supposons le contraire, \`a savoir qu'il existe des $\sigma_i' = \prod\limits_{\lambda\in\Lambda_i'} \lambda^{\varepsilon_\lambda} (\varepsilon_\lambda\in \{-1,0,1\})$ tels que $\prod \sigma_i'=1$. Soit $j$ un entier tel que $d(\sigma_i') |\Lambda_i|$ soit maximum pour $i=j$. Revenons \`a la d\'efinition de $\Lambda_j'$. Comme $(\sigma_j',(\sigma_k'))$ est un syst\`eme permis, on est dans le cas $\beta)$, o\`u l'on a d\'efini un syst\`eme permis $(\sigma,(\rho_k))$ tel que $d(\sigma)$ soit maximum. Posons $\sigma'=\sigma\sigma_j'$ et $\rho_k'=\rho_k\sigma_k'$ $(k\neq j)$. Comme les ``supports''  de $\sigma$ et de $\sigma_j'$ sont disjoints, on a $d(\sigma')=d(\sigma)+d(\sigma_j')$. Pour $k\neq j$, on a $d(\rho_k')\le d(\rho_k) +d(\sigma_k')$, et on sait que $d(\rho_k)|\Lambda_k| \le d(\sigma)|\Lambda_j|$ et $d(\sigma_k') |\Lambda_k| \le d(\sigma_j')|\Lambda_j|$, donc
$$
d(\rho_k')|\Lambda_k| \le d(\sigma)|\Lambda_j| + d(\sigma_j')|\Lambda_j| = d(\sigma')|\Lambda_j|\,.
$$
Le syst\`eme $(\sigma',(\rho_k'))$ est donc permis, et $d(\sigma') > d(\sigma)$, contrairement \`a la d\'efinition de $\sigma$. La contradiction \'etablit que $\bigcup \Lambda_i'$ est $qi$. Ainsi se termine la preuve par Bourgain que la propri\'et\'e de Sidon est \'equivalente \`a chacune des conditions \textbf{(CR), (CP)} et \textbf{(CB).}
\hfill$\blacksquare$

\vskip2mm

\hfill Jean--Pierre Kahane 05.08.2007

 \eject
 
 \vglue 2cm

\vskip4mm

\hfill\begin{minipage}{6,5cm}
Jean--Pierre Kahane

Laboratoire de Math\'ematique

Universit\'e Paris--Sud, B\^at. 425

91405 Orsay Cedex

\textsf{Jean-Pierre.Kahane@math.u-psud.fr}

\end{minipage}

\end{document}